\documentclass[12pt]{amsart}
\usepackage{graphicx}
\usepackage{amssymb}
\usepackage{amsfonts}
\usepackage{amsmath}
\usepackage{array}
\usepackage{rotating}

\newtheorem{example}{Example}[section]

\newtheorem{theorem}[example]{Theorem}

\newtheorem{proposition}[example]{Proposition}

\def\J{{\bf J}}
\def\FQSym{{\bf FQSym}}

\def\PQSym{{\bf PQSym}}
\def\CQSym{{\bf CQSym}}

\def\ev{{\rm Ev}}

\def\ssh{\Cup}

\def\sconc{\bullet}
\def\Std{{\rm Std}}
\def\Park{{\rm Park}}

\def\convol{{*}}

\def\park{{\bf a}}

\def\<{\langle}
\def\>{\rangle}

\def\Z{\operatorname{\mathbb Z}}

\def\F{{\bf F}}

\def\G{{\bf G}}

\def\P{{\bf P}}

\def\SG{{\mathfrak S}}

\def\Sym{{\bf Sym}}

\def\PF{{\rm PF}}
\def\PPF{{\rm PPF}}
\def\MM{{\mathcal M}}
\def\RR{{\bf R}}
\def\shuff#1#2{\mathbin{
\hbox{\vbox{ \hbox{\vrule \hskip#2 \vrule height#1 width 0pt
}%
\hrule}%
\vbox{ \hbox{\vrule \hskip#2 \vrule height#1 width 0pt
\vrule }%
\hrule}%
}}}

\def\shuf{{\mathchoice{\shuff{7pt}{3.5pt}}%
{\shuff{6pt}{3pt}}%
{\shuff{4pt}{2pt}}%
{\shuff{3pt}{1.5pt}}}}%
\def\shuffle{\,\shuf\,}

\def\Tabvrule{\vrule width-0.4pt}       % Difference de largeur
\def\Tabhrule{\hrule \hrule height-0.4pt} % Difference de hauteur
\def\Tabstrut{\vrule height2.2ex % Sur la ligne
                     depth0.8ex  % Sous la ligne
                     width0ex    % centrage horizontal
\relax}

\def\PasCase#1{\omit%
            $\vcenter{\hbox {\vbox to 0.4pt{}}
               \hbox{\makebox[3ex]{\Tabstrut$#1$}}}%
               \Tabvrule$}
\def\PasCasePoint{\PasCase{\cdot}}
\def\DessinCarre#1{%
    \vcenter{\hbox{}\hrule
             \hbox{\vrule\makebox[3ex]{\Tabstrut$#1$}\vrule}\Tabhrule}%
             \Tabvrule}
\def\GenRuban#1{\vcenter{\halign{&$\DessinCarre{##}$\cr#1}}\egroup}

\def\sTabvrule{\vrule width-0.4pt}
\def\sTabhrule{\hrule \hrule height-0.4pt}
\def\sTabstrut{\vrule height1.6ex depth0.6ex width0ex \relax}
\def\sDessinCarre#1{%
    \vcenter{\hbox{}\hrule
             \hbox{\vrule\makebox[2.3ex]%
                  {\sTabstrut$\scriptstyle#1$}\vrule}\sTabhrule}%
             \sTabvrule}
\def\sGenRuban#1{\vcenter{\halign{&$\sDessinCarre{##}$\cr#1}}\egroup}

\def\ruban{%
  \bgroup
  \let\ =\omit
  \let\\=\cr
  \let\x=\times
  \let\.=\PasCasePoint
  \offinterlineskip
  \GenRuban}

\def\sruban{%
  \bgroup
  \let\ =\omit
  \let\x=\times
  \let\\=\cr
  \offinterlineskip
  \sGenRuban}

\title{Parking Functions and Descent Algebras}

\author[J.-C.~Novelli and J.-Y.~Thibon]%
{Jean-Christophe Novelli and Jean-Yves Thibon}

\address[] {Institut Gaspard Monge, Universit\'e de Marne-la-Vall\'ee \\
5 Boulevard Descartes \\Champs-sur-Marne \\77454 Marne-la-Vall\'ee cedex 2 \\
FRANCE}
\email[Jean-Christophe Novelli]{novelli@univ-mlv.fr}
\email[Jean-Yves Thibon]{jyt@univ-mlv.fr} 
\date{}

\begin{document}

\begin{abstract}
We show that the notion of parkization of a word, a variant of the classical
standardization, allows to introduce an internal product on 
the Hopf algebra of parking functions. Its Catalan subalgebra is
stable under this operation and contains the descent algebra 
as a left ideal.
\end{abstract}

\maketitle

%%%%%%%%%%%%%%%%%%%%%%%%%%%%%%%%%%%%%%%%%%%%%%%%%%%%%%%%%%%%%%%%%%%%%%%%%%%%%%%
%%%%%%%%%%%%%%%%%%%%%%%%%%%%%%%%%%%%%%%%%%%%%%%%%%%%%%%%%%%%%%%%%%%%%%%%%%%%%%%
%%%%%%%%%%%%%%%%%%%%%%%%%%%%%%%%%%%%%%%%%%%%%%%%%%%%%%%%%%%%%%%%%%%%%%%%%%%%%%%
\section{Introduction}

Solomon \cite{Sol} constructed for each finite Coxeter group a remarkable
subalgebra of its group algebra, now called its descent algebra.

For the infinite series of Weyl groups, the direct sums of descent algebras
can be endowed with some interesting extra structure.  This is most
particularly the case for symmetric groups (type $A$), where the direct sum
$\Sigma=\bigoplus_{n\ge 0}\Sigma_n$ ($\Sigma_n$ being the descent algebra of
$\SG_n$) builds up a Hopf algebra, isomorphic to $\Sym$ (noncommutative
symmetric functions) and dual to $QSym$ (quasi-symmetric functions).

It has been understood by Reutenauer \cite{Re} and Patras \cite{Pa} that
$\Sigma$ could be interpreted as a subalgebra of the direct sum
$\SG=\sum_{n\ge 0}\Z\SG_n$ for the \emph{convolution product} of permutations,
which arises when permutations are regarded as graded endomorphisms of a free
associative algebra. Indeed, $\Sigma$ is then just the convolution subalgebra
generated by the homogeneous components of the identity map.
Further understanding of the situation has been provided by Malvenuto and
Reutenauer \cite{MR}, who gave a complete description of the Hopf algebra
structure of $\SG$, and by Poirier-Reutenauer \cite{PR}, who discovered an
interesting subalgebra based on standard Young tableaux.

Finally, the introduction of the Hopf algebra of free quasi-symmetric
functions $\FQSym$ \cite{NCSF6} clarified the picture and brought up a great
deal of simplification.  Indeed, $\FQSym$ is an algebra of noncommutative
polynomials over some auxiliary set of variables $a_i$, which is isomorphic to
$\SG$, and is mapped onto ordinary quasi-symmetric function $QSym$ when the
$a_i$ are specialized to commuting variables $x_i$, the natural basis
$\F_\sigma$ of $\FQSym$ going to Gessel's fundamental basis $F_I$. At the
level of $\FQSym$, the coproduct has a transparent definition (ordered
sum of alphabets), and most of its properties become obvious.

There is at least one point, however, on which this construction does not
shed much light. It is the original product of the descent algebras
$\Sigma_n$, which gives rise on $\Sym$ to a noncommutative analogue of the
internal product of symmetric functions (see~\cite{Mcd} for the classical
case). The introduction of the Hopf structure of $\Sigma=\Sym$ was extremely
useful, thanks to the so-called \emph{splitting formula} \cite{NCSF1,NCSF2}, a
compatibility property between all operations (internal and external product,
coproduct).
But the embedding of $\Sym$ in $\FQSym$ does not seem to bring new
information. In particular, the coproduct dual to the composition of
permutations has no nice definition in terms of product of alphabets, and the
splitting formula is no more valid in general. Hopf subalgebras in which it
remains valid have been studied by Schocker (Lie idempotent algebra,
\cite{Scho}) and by Patras-Reutenauer \cite{PR}, this last one being maximal
with respect to this property.

There are many combinatorial objects which can be regarded, in one way or
another, as generalizations of permutations.  Among them are \emph{parking
functions}, on which a Hopf algebra structure $\PQSym$, very similar to that
of $\FQSym$, can be defined \cite{NT1}. Actually, $\FQSym$ is a Hopf
subalgebra of $\PQSym$.

The aim of this note is to show that it is possible to define on $\PQSym$ an
internal product, dual to a natural coproduct corresponding to the Cartesian
product of ordered alphabets, exactly as in Gessel's construction of the
descent algebra \cite{Ges}.  This product is very different from the
composition  permutations or endofunctions, and looks actually rather strange.
It can be characterized in terms of the fundamental notion of
\emph{parkization} of words defined over a totally ordered alphabet in which
each element has a successor.

In \cite{NT1}, various Hopf subalgebras of $\PQSym$ have been introduced. We
shall show that the Catalan subalgebra $\CQSym$ (based on the Catalan family
of nondecreasing parking functions, or equivalently, non-crossing partitions)
is stable under this new internal product, and contains the descent algebra as
a left ideal. Moreover, the splitting formula remains valid for it.

\medskip
{\footnotesize
{\it Acknowledgements.-}
This project has been partially supported by EC's IHRP Programme, grant
HPRN-CT-2001-00272, ``Algebraic Combinatorics in Europe".
The authors would also like to thank the contributors of the
MuPAD project, and especially of the
combinat part, for providing the development environment for their research.
}

%%%%%%%%%%%%%%%%%%%%%%%%%%%%%%%%%%%%%%%%%%%%%%%%%%%%%%%%%%%%%%%%%%%%%%%%%%%%%%%
%%%%%%%%%%%%%%%%%%%%%%%%%%%%%%%%%%%%%%%%%%%%%%%%%%%%%%%%%%%%%%%%%%%%%%%%%%%%%%%
%%%%%%%%%%%%%%%%%%%%%%%%%%%%%%%%%%%%%%%%%%%%%%%%%%%%%%%%%%%%%%%%%%%%%%%%%%%%%%%
\section{Parking functions and parkization}

A \emph{parking function} on $[n]=\{1,2,\ldots,n\}$ is a word
$\park=a_1a_2\cdots a_n$ of length $n$ on $[n]$ whose non-decreasing
rearrangement $\park^\uparrow=a'_1a'_2\cdots a'_n$ satisfies $a'_i\le i$ for
all $i$.
Let $\PF_n$ be the set of such words.

One says that $\park$ has a \emph{breakpoint} at $b$ if $|\{\park_i\le b\}|=b$.
Then, $\park\in \PF_n$ is said to be \emph{prime} if its only breakpoint is
$b=n$.
Let $\PPF_n\subset\PF_n$ be the set of prime parking functions on $[n]$.
%It can easily be shown that $|\PPF_n|=(n-1)^{n-1}$ (see~\cite{Stan2}).

For a word $w$ on the alphabet ${1,2,\ldots}$, denote by $w[k]$ the word
obtained by replacing each letter $i$ by $i+k$.
If $u$ and $v$ are two words, with $u$ of length $k$, one defines
the \emph{shifted concatenation}
\begin{equation}
u\sconc v = u\cdot (v[k])
\end{equation}
and the \emph{shifted shuffle}
\begin{equation}
u\ssh v= u\shuffle (v[k])\,.
\end{equation}
The set of permutations is closed under both operations, and the subalgebra
spanned by this set is isomorphic to $\SG$ \cite{MR} or to $\FQSym$
\cite{NCSF6}.

Clearly, the set of all parking functions is also closed under these
operations.  The prime parking functions exactly are those which do not occur
in any nontrivial shifted shuffle of parking functions. These properties
allowed us to define a Hopf algebra of parking functions in \cite{NT1}.

This algebra, denoted by $\PQSym$, for \emph{Parking Quasi-Symmetric
functions}, is spanned as a vector space by elements $\F_\park$
($\park\in\PF$), the product being defined by
\begin{equation}
\label{prodF}
\F_{\park'}\F_{\park''}:=\sum_{\park\in\park'\ssh\park''}\F_\park\,.
\end{equation}

For example,
\begin{equation}
\F_{12}\F_{11}= \F_{1233} + \F_{1323} + \F_{1332} + \F_{3123} + \F_{3132}
+ \F_{3312}\,.
\end{equation}

The coproduct on $\PQSym$ is a natural extension of that of $\FQSym$.
Recall (see~\cite{MR,NCSF6}) that if $\sigma$ is a permutation,
\begin{equation}
\label{CoprodF}
\Delta\F_{\sigma} = \sum_{u\cdot v=\sigma}{\F_{\Std(u)} \otimes \F_{\Std(v)}},
\end{equation}
where $\Std$ denotes the usual notion of standardization of a word.

For a word $w$ over a totally ordered alphabet in which each element
has a successor, we defined in \cite{NT1} a
notion of \emph{parkized word} $\Park(w)$, a parking function which reduces to
$\Std(w)$ when $w$ is a word without repetition.

For $w=w_1w_2\cdots w_n$ on $\{1,2,\ldots\}$, we set
\begin{equation}
\label{dw}
d(w):=\min \{i | \#\{w_j\leq i\}<i \}\,.
\end{equation}
If $d(w)=n+1$, then $w$ is a parking function and the algorithm terminates,
returning~$w$. Otherwise, let $w'$ be the word obtained by decrementing all
the elements of $w$ greater than $d(w)$. Then $\Park(w):=\Park(w')$. Since
$w'$ is smaller than $w$ in the lexicographic order, the algorithm terminates
and always returns a parking function.

\smallskip
For example, let $w=(3,5,1,1,11,8,8,2)$. Then $d(w)=6$ and the word
$w'=(3,5,1,1,10,7,7,2)$.
Then $d(w')=6$ and $w''=(3,5,1,1,9,6,6,2)$. Finally, $d(w'')=8$ and
$w'''= (3,5,1,1,8,6,6,2)$, that is a parking function.
Thus, $\Park(w)=(3,5,1,1,8,6,6,2)$.

\smallskip
The coproduct on $\PQSym$ is defined by
\begin{equation}
\Delta \F_{\park}:= \sum_{u\cdot v=\park} \F_{\Park(u)} \otimes \F_{\Park(v)},
\end{equation}
For example,
\begin{equation}
\Delta\F_{3132} = 1\otimes\F_{3132} + \F_{1}\otimes\F_{132} +
\F_{21}\otimes\F_{21} + \F_{212}\otimes\F_{1} + \F_{3132}\otimes 1\,.
\end{equation}
The product and the coproduct of $\PQSym$ are compatible, so
that $\PQSym$ is a graded bialgebra, connected, hence a Hopf algebra.
Let $\G_{\park}=\F_{\park}^* \in\PQSym^*$ be the dual basis of $(\F_\park)$.
If $\langle\,,\,\rangle$ denotes the duality bracket, the product on
$\PQSym^*$ is given by
\begin{equation}
\label{prodG}
\G_{\park'} \G_{\park''} = \sum_{\park}
    \langle\, \G_{\park'}\otimes\G_{\park''}, \Delta\F_\park \,\rangle\,
    \G_\park
= \sum_{\park \in \park'\convol\park''} \G_\park\,,
\end{equation}
where the \emph{convolution} $\park'\convol\park''$ of two parking functions
is defined as
\begin{equation}
\park'\convol\park'' = \sum_{u,v ;
\park=u\cdot v, \Park(u)=\park', \Park(v)=\park''} \park\,.
\end{equation}
For example,
\begin{equation}
\begin{split}
\G_{12} \G_{11} &= \G_{1211} + \G_{1222} + \G_{1233} + \G_{1311} + \G_{1322}\\
&+ \G_{1411} + \G_{1422} + \G_{2311} + \G_{2411} + \G_{3411}\,.
\end{split}
\end{equation}
When restricted to permutations, it coincides with the convolution
of~\cite{Re,MR}.

The coproduct of a $\G_\park$ is
\begin{equation}
\Delta \G_\park := \sum_{u,v ; \park\in u\ssh v}
                   {\G_{\Park(u)} \otimes \G_{\Park(v)}}\,.
\end{equation}
For example,
\begin{equation}
\begin{split}
\Delta \G_{41252} &= 1 \otimes\G_{41252} + \G_{1}\otimes\G_{3141} +
                   \G_{122}\otimes\G_{12} \\
                 &+ \G_{4122}\otimes\G_{1} + \G_{41252}\otimes1\,.
\end{split}
\end{equation}

%%%%%%%%%%%%%%%%%%%%%%%%%%%%%%%%%%%%%%%%%%%%%%%%%%%%%%%%%%%%%%%%%%%%%%%%%%%%%%%
%%%%%%%%%%%%%%%%%%%%%%%%%%%%%%%%%%%%%%%%%%%%%%%%%%%%%%%%%%%%%%%%%%%%%%%%%%%%%%%
%%%%%%%%%%%%%%%%%%%%%%%%%%%%%%%%%%%%%%%%%%%%%%%%%%%%%%%%%%%%%%%%%%%%%%%%%%%%%%%
\section{Polynomial realization of $\PQSym^*$}

In the sequel, we need the following definitions: given a totally ordered
alphabet $A$, the \emph{evaluation vector} $Ev(w)$ of a word $w$ is the
sequence of number of occurrences of all the elements of $A$ in $w$. The
\emph{packed evaluation vector} $c(w)$ of $w$ is obtained from $Ev(w)$ by
removing all its zeros. The \emph{fully unpacked evaluation vector} $d(w)$ of
$w$ is obtained from $c(w)$ by inserting $i-1$ zeros after each entry $i$ of
$c(w)$ except the last one.
For example, if $w=3117291781329$, $Ev(w)=(4,2,2,0,0,0,2,1,2)$,
$c(w)=(4,2,2,2,1,2)$, and $d(w)=(4,0,0,0,2,0,2,0,2,0,1,2)$.

The algebra $\PQSym^*$ admits a simple realization in terms of noncommutative
polynomials, which is reminescent of the construction of $\FQSym$. 
If $A$ is a totally ordered infinite alphabet, one can set
\begin{equation}\label{realG}
\G_\park(A)=\sum_{w\in A^*, \Park(w)=\park}w
\end{equation}
These polynomials satisfy the relations (\ref{prodG}) and allow to write the
coproduct as $\Delta \G_a = \G_a(A'\hat{+} A'')$ where $A'\hat{+} A''$ denotes
the ordered sum of two mutually commuting alphabets isomorphic to $A$ as
ordered sets.

Recall from \cite{NT1} that the sums
\begin{equation}
\P^\pi := \sum_{\park ; \park^\uparrow=\pi} {\F_\park}
\end{equation}
where $\park^\uparrow$ means the non-decreasing reordering and $\pi$ runs over
non-decreasing parking functions, span a cocommutative Hopf subalgebra
$\CQSym$ of $\PQSym$.

As with $\FQSym$, one can take the commutative image of the $\G_\park$, that
is, replace the alphabet $A$ by an alphabet $X$ of commuting variables
(endowed with an isomorphic ordering). Then, $\G_{\park'}(X)=\G_{\park''}(X)$
iff $\park'$ and $\park''$ have the same non-decreasing reordering $\pi$, and
both coincide with the generalized quasi-monomial function
$\MM_\pi=(\P^\pi)^*$ of \cite{NT1}, that is, the natural basis of the
commutative Catalan algebra $\CQSym^*$.

Actually, $\CQSym^*$ contains $QSym$ as a subalgebra, the quasi-monomial
functions being obtained as $M_I=\sum_{c(\pi)=I}\MM_\pi$.
%, where the sum runs
%over
%all non-decreasing parking functions which evaluation vector after discarding
%their zeroes is $I$. 

As a first application of the polynomial realization, we can quantize
$\CQSym^*$. Indeed, we can proceed as for the quantization of $QSym$
\cite{TU}, that is, we map the $a_i$ on $q$-commuting variables $x_i$, that
is, $x_jx_i=q x_i x_j$ for $i<j$, $\G_{\park'}(X)$ and $\G_{\park''}(X)$ are
equal only up to a power of $q$ when $\park'$ and $\park''$ have the same
non-decreasing reordering $\pi$, and the resulting algebra is not commutative
anymore. Deforming the coproduct so as to maintain compatibility with the
product, we obtain a self-dual Hopf algebra, which is isomorphic to the
Loday-Ronco algebra of plane binary trees \cite{LR}.

However, our main application will be the definition of an internal product on
$\PQSym$.

%%%%%%%%%%%%%%%%%%%%%%%%%%%%%%%%%%%%%%%%%%%%%%%%%%%%%%%%%%%%%%%%%%%%%%%%%%%%%%%
%%%%%%%%%%%%%%%%%%%%%%%%%%%%%%%%%%%%%%%%%%%%%%%%%%%%%%%%%%%%%%%%%%%%%%%%%%%%%%%
%%%%%%%%%%%%%%%%%%%%%%%%%%%%%%%%%%%%%%%%%%%%%%%%%%%%%%%%%%%%%%%%%%%%%%%%%%%%%%%
\section{The internal product}

Recall that Gessel constructed the descent algebra by extending to $QSym$ the
coproduct dual to the internal product of symmetric functions. That is, if $X$
and $Y$ are two totally and isomorphically ordered alphabets of commuting
variables, we can identify a tensor product $f\otimes g$ of quasi-symmetric
functions with $f(X)g(Y)$. Denoting by $XY$ the Cartesian product $X\times Y$
endowed with the lexicographic order, Gessel defined for $f\in QSym_n$
\begin{equation}
\delta(f)=f(XY) \in QSym_n\otimes QSym_n \,.
\end{equation}
The dual operation on $\Sym_n$ is the internal product $*$, for which $\Sym_n$
is anti-isomorphic to the descent algebra $\Sigma_n$.

This construction can be extended to the commutative Catalan algebra
$CQSym^*$, and in fact, even to $\PQSym^*$.

Let $A'$ and $A''$ be two totally and isomorphically ordered alphabets of 
noncommuting variables, but such that $A'$ and $A''$ commute with each other.
We denote by $A'A''$ the Cartesian product $A'\times A''$ endowed with the
lexicographic order. This is a total order in which
each element has a successor, so that $G_\park(A'A'')$ is a well
defined polynomial. Identifying tensor products $u\otimes v$ of words of the
same length with words over $A'A''$, we have
\begin{equation}
G_\park(A'A'')=\sum_{\Park(u\otimes v)=\park} u\otimes v\,.
\end{equation}
Our main result is the following
\begin{theorem}
The formula $\delta(G_\park)=G_\park(A'A'')$ defines a coassociative coproduct
on each homogeneous component $\PQSym_n^*$. Actually,
\begin{equation}
\delta(\G_\park)=\sum_{\Park(\park'\otimes\park'')=\park}
                     {\G_{\park'}\otimes \G_{\park''}}\,.
\end{equation}
By duality, the formula
\begin{equation}
\F_{\park'} * \F_{\park''} = \F_{\Park(\park'\otimes\park'')}
\end{equation}
defines an associative product on each $\PQSym_n$.
\end{theorem}

\begin{example}%{\rm }
\begin{equation}
\F_{211}*\F_{211}=\F_{311}; \qquad \F_{211}*\F_{112} = \F_{312};
\end{equation}
\begin{equation}
\F_{211}*\F_{121}=\F_{321}; \qquad \F_{112}*\F_{312} = \F_{213};
\end{equation}
\begin{equation}
\F_{31143231}*\F_{23571713} = \F_{61385451}.
\end{equation}
\end{example}

%{\footnotesize Formule d'\'eclatement vraie \`a ce niveau ?} -> non
 
%%%%%%%%%%%%%%%%%%%%%%%%%%%%%%%%%%%%%%%%%%%%%%%%%%%%%%%%%%%%%%%%%%%%%%%%%%%%%%%
%%%%%%%%%%%%%%%%%%%%%%%%%%%%%%%%%%%%%%%%%%%%%%%%%%%%%%%%%%%%%%%%%%%%%%%%%%%%%%%
%%%%%%%%%%%%%%%%%%%%%%%%%%%%%%%%%%%%%%%%%%%%%%%%%%%%%%%%%%%%%%%%%%%%%%%%%%%%%%%
\section{Subalgebras of $(\PQSym_n,*)$}

The following result is almost immediate.

\begin{proposition}
The homogeneous components $\CQSym_n$ of the Catalan algebra are
stable under the internal product $*$.
\end{proposition}

\begin{example}
\begin{equation}
\P_{1123} * \P_{1111} = \P_{1134}; \qquad
\P_{1111} * \P_{1123} = \P_{1123}.
\end{equation}
\begin{equation}
\P_{1123} * \P_{1112} = 2\P_{1134} + \P_{1234}; \qquad
\P_{1122} * \P_{1224} = \P_{1134} + \P_{1233} + 2\P_{1234}.
\end{equation}
\begin{equation}
\P_{1123} * \P_{1224} = 2\P_{1134} + 5\P_{1234}.
\end{equation}
%\begin{equation}
%\end{equation}
\end{example}

It is interesting to observe that these algebras are non-unital.
Indeed, as one can see on the first example just above

\begin{proposition}
The element $\J_n=\P^{(1^n)}$ is a left unit for $*$, but not
a right unit.
\end{proposition}

%\begin{equation}
%\end{equation}
The splitting formula is valid in $\CQSym_n$. That is,

\begin{proposition}
Let $\mu_r$ denote the $r$-fold product map from $\CQSym^{\otimes r}$
to $\CQSym$, $\Delta^r$ the $r$-fold coproduct with values in
$\CQSym^{\otimes r}$, and $*_r$ the internal product of the $r$-fold tensor
product of algebras $\CQSym^{\otimes r}$. Then, for
$f_1,\ldots,f_r,g\in\CQSym$,
\begin{equation}
(f_1\cdots f_r)*g=\mu_r[(f_1\otimes\cdots\otimes f_r)*_r \Delta^r(g)]\,.
\end{equation}
\end{proposition}
This is exactly the same formula as with the internal product of $\Sym$,
actually, an extension of it, since we have

\begin{theorem}
The Hopf subalgebra of $\CQSym$ generated by the elements $\J_n$, which is
isomorphic to $\Sym$ by $j:S_n\mapsto \J_n$, is stable under $*$, and thus
also $*$-isomophic to $\Sym$.
Moreover, the map $f\mapsto f*\J_n$ is a projector onto $\Sym_n$, which is
therefore a left $*$-ideal of $\CQSym_n$.
\end{theorem}

If $i<j<\ldots<r$ are the letters occuring in $\pi$, so that as a word
$\pi=i^{m_i}j^{m_j}\cdots r^{m_r}$, then
\begin{equation}
\P^\pi*\J_n = \J_{m_i}\J_{m_j}\cdots \J_{m_r}\,.
\end{equation} 

In the classical case, the non-commutative complete fonctions split into a sum
of ribbon Schur functions, using a simple order on compositions.
To get an analogous construction in our case, we have defined a partial order
on non-decreasing parking functions.

Let $\pi$ be a non-decreasing parking function and $\ev(\pi)$ be its
evaluation vector. The successors of $\pi$ are the non-decreasing parking
functions whose evaluations are given by the following algorithm: given two
non-zero elements of $\ev(\pi)$ with only zeros between them, replace the
left one by the sum of both and the right one by 0.
For example, the successors of $113346$ are $111146$, $113336$, and $113344$.

By transitive closure, the successor map gives rise to a partial order on
non-decreasing parking functions. We will write $\pi\preceq\pi'$ if $\pi'$ is
obtained from $\pi$ by successive applications of successor maps.

The Catalan ribbon functions are defined by
\begin{equation}
\label{catalRub}
\P^\pi =: \sum_{\pi'\succeq\pi} {\RR_{\pi'}}\,.
\end{equation}
%This characterizes the $\RR_{\pi}$.
%
%\smallskip
%The product of two $\RR$ functions is then
%\begin{equation}
%\RR_{\pi'} \RR_{\pi''} = \RR_{\pi'\sconc\pi''} +
%\RR_{\pi'\triangleright\pi''}\,,
%\end{equation}
%where $\triangleright$ denotes the shifted concatenation defined by shifting
%all elements of $\pi''$ by the difference between the greatest and the
%smallest element of $\pi'$.
%
%For example,
%\begin{equation}
%\RR_{11224} \RR_{113} = \RR_{11224668} + \RR_{11224446}\,.
%\end{equation}

The $\RR_{\pi}$ are the pre-images of the ordinary ribbons under the
projection $f\mapsto f*\J_n$:
\begin{proposition}
Let $I$ be the composition obtained by discarding the zeros of the evaluation
of an non-decreasing parking function $\pi$. Then
\begin{equation}
\RR_\pi * \J_n = j(R_I).
\end{equation}
More precisely, if $I=(i_1,\ldots,i_p)$, this last element is equal to
$\RR_{1^{i_1}\sconc1^{i_2}\sconc\cdots\sconc1^{i_p}}$, that is, the Catalan
ribbon indexed by the only non-decreasing word of evaluation $d(\pi)$.
\end{proposition}

%%%%%%%%%%%%%%%%%%%%%%%%%%%%%%%%%%%%%%%%%%%%%%%%%%%%%%%%%%%%%%%%%%%%%%%%%%%%%%%
%%%%%%%%%%%%%%%%%%%%%%%%%%%%%%%%%%%%%%%%%%%%%%%%%%%%%%%%%%%%%%%%%%%%%%%%%%%%%%%
%%%%%%%%%%%%%%%%%%%%%%%%%%%%%%%%%%%%%%%%%%%%%%%%%%%%%%%%%%%%%%%%%%%%%%%%%%%%%%%
%%%%%%%%%%%%%%%%%%%%%%%%%%%%%%%%%%%%%%%%%%%%%%%%%%%%%%%%%%%%%%%%%%%%%%%%%%%%%%%

\end{document}